\documentclass[a4paper, 12pt, reqno]{amsart}
\usepackage{amsmath, amssymb, caption, latexsym, graphicx}
\usepackage[ansinew]{inputenc}
\usepackage{fontenc}
\usepackage[dvips, dvipsnames, usenames]{color}
\usepackage[stable]{footmisc}
\usepackage[backref]{hyperref}
\hypersetup{colorlinks=true, urlcolor=blue, citecolor=blue}

\newtheorem{theorem}{Theorem}

\begin{document}

\title{Beyond Chaos}

\maketitle

\begin{center}

\large{Introducing recursive logistic interactions}

\vspace{20pt}

\small{

Antonio Leon Sanchez\\
I.E.S Francisco Salinas, Salamanca, Spain\\
\href{http://www.interciencia.es}{http://www.interciencia.es}\\
\href{mailto:aleon@interciencia.es}{aleon@interciencia.es}\\
}
\end{center}


\begin{abstract}
The first part of this paper defines recursive interactions by means of logistic
functions and derives a general result on the way interacting systems evolve in
attractors. It also defines the notion of coevolution trajectory and presents a new
family of attractors: orbital attractors (including single, irregular, folded, complex
and discontinuous orbits). The second part summarizes the results of a first
experimental analysis of recursive interactions in both binary and multiple
interactions. Among other results, this analysis reveals that interacting systems may
easily evolve from chaos to order.
\end{abstract}


\section{Introduction}

\noindent Mathematical continuity makes no sense in biology. The biosphere is
radically discontinuous, i.e. discrete; it is the paradigm of discreteness. Each
living being is a discrete unity, unique and unrepeatable, that emerges and self
maintains at the expense of a discrete net of metabolic reactions governed by a
discrete net of discrete information units. Each living being forms part of a discrete
group of discrete individuals which in turns form part of a discrete net of ecological
interactions. But living organisms are not the only discrete objects in the universe,
matter and energy are also discrete. In consequence, all natural processes involving
matter or energy interchanges have to be of a discrete nature. Even space and time
could be discrete as has been repeatedly suggested from different areas of physics
(\cite{Biser1941}, \cite{Coish1959}, \cite{Finkelstein1986}, \cite{Meessen1989},
\cite{Kragh1994}, \cite{Forrest1995}, \cite{Susskind1997}, \cite{Bendegem1997},
\cite{Susskind1997}, \cite{Bekenstein1980}, \cite{Rovelli2001}, \cite{Green2001},
\cite{Bendegem2002}, \cite{Bekenstein2003}, \cite{Smolin2003}, \cite{Baez2003} ,
\cite{Green2004}, \cite{Veneziano2004}, \cite{Smolin2004}, \cite{Lloyd2005},
\cite{Bekenstein2003}). In these conditions, and being the best model of any object
the very object itself, what should be put into question is not the role of discrete
models \cite{Deng2007} but the role of mathematical continuity in an essentially
discontinuous world.

In addition of being discrete and interactive, living organisms are also recursive.
Each generation determines the next one through a reproductive process that is not
exact, giving therefore the appropriate opportunity to evolution. It is then clair the
reason for which biology has always been interested in discrete and recursive models.
Among them the logistic function, perhaps the most simple and productive model of some
biological significance. From the pioneering works of Stanislaw Ulam, Paul Stein
\cite{Ulam1963}, Nicholas C. Metropolis \cite{Metropolis1973} and Robert May
\cite{May1976} to nowadays, iterative calculus and the logistic function have been an
inexhaustible source of surprising results, from deterministic chaos to periodic
attractors; results that, on the other hand, were immediately generalized by the so
called \emph{principle of universality}. As is well known, the control parameter
$\lambda$ of a logistic function:
\begin{equation}
x_{n+1} = 4 \lambda x_n (1 - x_n)
\end{equation}

\noindent determines the system evolution, particularly the type of attractor it
finally falls in. Now then, biological systems do not evolve separately but immerse in
a complex interaction network. For this reason biology has also been interested in
modeling interactions. Most of those models, as the classical Lotka-Volterra, are
defined in terms of differential (or finite differences) equations. But as far as I
know, biological interactions have never been modeled by logistic functions. As we
will immediately see, biological (and non biological) interactions can easily be
modeled by means of this type of functions. And the results cannot be more interesting
from both the mathematical and the biological point of view. In fact, on the one hand
recursive logistic interactions are a new source of new mathematical objects as
\emph{coevolution trajectories} or orbital attractors. On the other, they provide us a
new way of examining coevolution processes and derive some significant results on the
way systems coevolve. For instance, that interacting systems may coevolve from chaotic
regimes to stable states defined by periodic attractors of low period; or that complex
attractors, as fractal or chaotic, may suddenly evolve in single periodic attractors;
or that systems to the brink of extinction may be completely recovered and stabilized
thanks to its recursive interactions with other systems; or, in the case of multiple
interactions, that system stability grows with the number of interactions.

The discussion that follows begins by defining the appropriate theoretical notions,
from which a general result related to the way interacting systems evolve in
attractors is formally derived. In the following two sections I define some new
mathematical objects as the coevolution trajectory of two interacting systems and a
new family of attractors. Finally, I resume the results of a first experimental work
on recursive logistic interactions, including a short presentation of multiple
interactions. Evidently, the field to explore is immense and what follows is but a
simple introductory note. All numerical and graphical data have been obtained with the
aid of
\href{http://www.interciencia.es/InterCalculus/InterCalculus.htm}{InterCalculus}
\cite{Leon2008b}, a computer application developed to analyze recursive interactions
by means of logistic functions.

\section{Recursive logistic interactions: Definitions}

\noindent Let us consider two systems $X$ and $Y$ modeled by two logistic real
functions:
\begin{align*}
&x_{n+1} = \lambda_x x_n (1 - x_n)\\
&y_{n+1} = \lambda_y y_n (1 - y_n)
\end{align*}

\noindent An elemental way to make both functions mutually dependent consists in
defining the control parameter of each function in terms of the other function
variable:
\begin{align*}
&x_{n+1} = g_{x}(y_n) x_n (1 - x_n)\\
&y_{n+1} = g_{y}(x_n) y_n (1 - y_n)
\end{align*}

\noindent If $g_{x}(y_n)$, the \emph{interaction function} of system $X$ with system
$Y$, is decreasing the interaction of system $X$ with system $Y$ will be said
negative; if it increasing the interaction will be said positive. The same applies to
$g_y(x_n)$, the \emph{interaction function} of system $Y$ with system $X$. There are
many ways to define the interaction functions $g_{x}$ and $g_{y}$. Here we will use
the following:
\begin{equation*}
g_{x}(y_n) = 4(1 - S_{x}y_n)
\end{equation*}

\noindent for negative interactions, and:
\begin{equation*}
g_{x}(y_n) = 4S_{x}y_n
\end{equation*}

\noindent for positive ones. In both cases, $S_{x}$ represents the sensitivity of
system $X$ to system $Y$, i.e. a measure of the effects of system $Y$ on system $X$.
It will take its values within the real half closed interval $(0, 1]$. Similarly, for
the interaction of system $Y$ with system $X$ we will have:
\begin{align*}
&g_{y}(x_n) = 4(1 - S_{y}x_n)\\
&g_{y}(x_n) = 4S_{y}x_n
\end{align*}

\noindent It is now immediate to define the following three types of recursive
logistic interactions (recursive interactions from now on):
\begin{enumerate}
\item Negative-Negative:
\begin{align}
&x_{n+1} = 4(1 - S_{x}y_n) x_n (1 - x_n) \label{eqn:neg-neg 1} \\
&y_{n+1} = 4(1 - S_{y}x_n) y_n (1 - y_n) \label{eqn:neg-neg 2}
\end{align}

\item Positive-Negative:
\begin{align}
&x_{n+1} = 4S_{x}y_n x_n (1 - x_n) \label{eqn:neg-pos 1}\\
&y_{n+1} = 4(1 - S_{y}x_n) y_n (1 - y_n) \label{eqn:neg-pos 2}
\end{align}

\item Positive-Positive:
\begin{align}
&x_{n+1} = 4 S_{x}y_n x_n (1 - x_n)\\
&y_{n+1} = 4 S_{y}x_n y_n (1 - y_n)
\end{align}
\end{enumerate}

\noindent Here we will exclusively deal with these recursive interactions. From now on
negative-negative interactions will be referred to as NN interactions;
positive-negative as PN; and positive-positive as PP.

The above definitions can immediately be generalized for any number of systems. In
effect, consider $m$ systems $\langle X_i \rangle_{i=1,2, \dots m}$ being each modeled
by a logistic function of a real variable $x_i$. Assume that all systems interact
which each other. Their evolution can be expressed as:
\begin{equation}\label{eqn:2-interaccion m-aria}
x_{i, \, n+1} = \frac{1}{m-1}\sum_{\underset{j\neq i}{j=1}}^{m} g_{x_i}(x_j) x_{i,n}
(1 - x_{i,n})
\end{equation}

\noindent where $x_{i,n}$ represent the $n$-th generation of system $x_i$;
$g_{x_i}(x_j)$ is the interaction function of system $X_i$ with system $X_j$, a
function that will be given by $4S_{i,j}x_{j,n}$ if that interaction is positive, or
$4(1 - S_{i,j}x_{j,n})$ if it is negative, being $S_{i,j}$ the sensitivity of system
$X_i$ to system $X_j$. Equations (\ref{eqn:2-interaccion m-aria}) can be easily
computerized so that it is possible to make both numerical and graphical analysis of
the evolution of any number of interacting systems.

\section{Attractors theorem}

\noindent We will prove now a general and basic result that states the simultaneity
and equivalence of the attractors the interacting systems evolve in.

\begin{theorem}
Let $X$ and $Y$ be two systems that interact with each other. A necessary and
sufficient condition for system $X$ to fall in a period-k attractor is that system $Y$
also falls in a period-k attractor.
\end{theorem}

\begin{proof}

Assume that system $X$ falls in a period-k attractor. From a certain integer $i > 0$
it will hold:
\begin{equation}\label{eqn:x_n = x_n+k}
x_n = x_{n+k}, \ \forall n > i
\end{equation}

\noindent On the other hand we will have
\begin{align*}
x_n = S(y_{n-1}) x_{n-1} (1 - x_{n-1})\\
x_{n+k}= S(y_{n+k-1}) x_{n+k-1} (1 - x_{n+k-1})
\end{align*}

\noindent where $S(y_{n-1})$ is $4(1 - S_xy_{n-1})$ if the interaction of $X$ with $Y$
is negative, or $4S_xy_{n-1}$ if it is positive. The same applies to $S(y_{n+k-1})$.
According to (\ref{eqn:x_n = x_n+k}) we can write:
\begin{equation*}
S(y_{n-1}) x_{n-1} (1 - x_{n-1}) = S(y_{n+k-1}) x_{n+k-1} (1 -
x_{n+k-1})
\end{equation*}

\noindent and taking into account that $x_{n-1} = x_{n+k-1}$ we have:
\begin{equation*}
S(y_{n-1}) = S(y_{n+k-1})
\end{equation*}

\noindent That is to say:
\begin{equation*}
4(1 - S_xy_{n-1}) = 4(1 - S_xy_{k+n-1})
\end{equation*}

\noindent if the interaction of $X$ with $Y$ is negative; or:
\begin{equation*}
4S_xy_{n-1} = 4S_xy_{k+n-1}
\end{equation*}

\noindent if it is positive. In both cases, we immediately get:
\begin{equation*}
y_{n-1} = y_{n+k-1}, \ \forall n > i
\end{equation*}

\noindent which implies that system $Y$ has also fallen in a period-k attractor. This
proves that if system $X$ evolves in a period-k attractor, so evolves system $Y$. A
similar argument proves the complementary result.

\end{proof}

\noindent The above theorem applies only to binary interactions. In the case of
multiple interactions, in fact, systems do not reach simultaneously their final
attractor; although, in most cases, once a first system falls in its final attractor
the remainder ones also fall in their respective final attractors after a few number
of iterations. I think attractors theorem is not the only conclusion that can be
formally derived for binary interactions. Experimental research strongly suggests
that, among others, the following results could also be formally proved:

\begin{itemize}

\item In NN interactions, the less sensitive system always goes ahead\footnote{System X
goes ahead of system Y if x > y.} of the most sensitive one.

\item In PN interactions, the system suffering the negative interaction always goes
ahead of the other.

\item Coevolution trajectories in NN interactions are symmetrical with respect the
bisectors of $(0, 1) \times (0, 1)$.

\end{itemize}

\section{Coevolution trajectories}

\noindent Recursive interaction invariably leads to a final attractor that may be
periodic, chaotic, fractal or of some new types, as we will immediately see. The set
of points each interacting systems traverses towards its final attractor defines its
coevolution trajectory (CT). It usually consists of two or more well defined lines,
branches from now on, usually with diagonal symmetry although they can also exhibit
radial or spiral symmetry. Each interacting system has its own CT, being both CT very
similar.

Systems evolve along the branches of their CT by cyclically stepping from a branch to
other, always in the same order, so that they exhibit a remarkable periodic behaviour.
Due to the discrete nature of iterations, CT are not continues but discrete lines,
although the density of points can be extremely high and we have to magnify them
thousands of millions times to make discontinuity visible. Naturally, in these
superdense regions the progress towards the final attractor is also extremely slow.
Coevolution trajectories are then sets of points -usually lines- in the interval $[0,
1] \times [0, 1]$ of $\mathbb{R}^2$ along which interacting systems evolve towards
their respective final attractors. As far as they have been examined, the coevolution
trajectories of both systems are geometrically similar, with the same number of
branches and the same type of symmetry. The branches of a CT may converge or diverge
from the central point $(0.5, 0.5)$, being the changes simultaneous in all branches of
both CT.

There also exists CT trough which systems evolve from chaos to non chaotic attractor
(periodic, fractal or orbital) or even from a chaotic regime to other different
chaotic regime. While traversing their respective CT, one of the systems goes always
ahead of the other (the less sensitive in the case of NN interactions, and the one
suffering the negative interaction in PN interactions). With the appropriate initial
seeds it is possible, however, to start from an inverted position. In these cases an
inversion of the CT occurs thanks to which the systems recover its normal relative
position. To go ahead of the other seems to be a general law that operates even in
chaotic regimes.

\section{Orbital attractors}

\noindent As in the case of single logistic recursion, recursive interactions also end
up by reaching a final attractor. Apart from the well known periodic, fractal and
chaotic attractors, in the case of recursive interactions we can also observe at least
a new family of attractors that we will term \emph{orbital attractors}. An orbital
attractor is typically composed of one or more closed lines (orbits) which are usually
eccentric. Each orbit its initially a discontinuous set of segments which
progressively extend in the same direction so that finally they overlap and close the
line. In some cases the orbits are single lines while in other they have a complex
internal structure. Systems trapped in an orbital attractor behave as if they were
orbiting, although their true behaviour is a little more complex. In fact, they jump
cyclically from an orbit to other always in the same order and in such a way that the
successive jumps on the same orbit take place on different segments, always in the
same order. All segments grow in the same direction so that after a few interactions
each segment reaches the next one, although, being of a discontinuous nature, their
respective points do not coincide, each segment continues indefinitely by occupying
the empty gaps along the orbit. After a considerable number of iterations, each
segment reaches itself and progress trough its own gaps. This way of progressing
suggest the possibility that after a huge number of iterations all orbital attractors
finally end in a periodic attractor.

The orbits of an orbital attractor may be at least of one of the following types:
\begin{enumerate}
\item{Single}
\item{Irregular}
\item{Folded}
\item{Complex}
\item{Discontinuous}
\end{enumerate}

\noindent The lines of a single orbit do not exhibit internal structure, they always
appear as single lines even if we magnify them by several billions times. The same
apply to irregular and discontinuous orbits although in these cases the lines exhibit
more or less irregular forms, as complex indentations. Folded orbit exhibit a complex
non fractal internal structure of folded lines in the regions of maximum curvature,
this structure seems disappear in the less curvature regions. This is also the case of
complex orbits, although in these orbits the internal structure is extremely complex,
perhaps of a fractal nature. In the case of discontinuous orbits, the segments do not
extend and maintain the empty gaps between them. After a certain number of iterations
these attractors degenerates in a single period attractor.

A remarkable characteristic of attractors in recursive logistic interactions is that
chaotic, fractal and orbital attractors may suddenly evolve in a periodic attractor
whose points presumably belong to the original attractor. It is also remarkable the
existence (at least in NN interactions) of periodic attractors whose attraction
decreases exponentially as systems approach them, so that it takes thousands of
millions of iterations to progress one decimal cipher towards the attractor final
value. They could be termed asymptotic attractors. It is possible that no finite
number of iterations suffices for the system to attaint the attractor (in the same
sense that the limit of a sequence cannot be reached by the successive terms of the
sequence).

\section{Negative-negative interactions}

\noindent In accordance with the above definition (\ref{eqn:neg-neg
1})-(\ref{eqn:neg-neg 2}), in NN interactions the growth of a systems is always to the
detriment of the other. Thus, NN interactions model the coevolution of two systems
that compete with each other. The results will depend on the sensitivities and on the
initial seeds of both systems, although the sensitivity dependence is stronger. The
coevolution of interacting systems is therefor controlled by four real number in the
interval $(0, 1]$, which means a huge number of possibilities to examine, each
representing a possible coevolution history. Although only an insignificant number of
cases has been examined, the following conclusions could be of general application:

\begin{enumerate}

\item The branches of the coevolution trajectories are symmetrical with respect to the
bisector, $x_{n+1} = x_{n}$ in the case of system $X$, or  $y_{n+1} = y_{n}$ in the
case of system $Y$).

\item All branches of a CT are geometrically similar, with same length and the same
point density.

\item Systems traverse the branches of their respective CT by stepping from one branch
to other, always in the same order so that they exhibit a remarkable periodic
behaviour while evolving towards the final attractor.

\item The oscillations in NN interactions are synchronic, i.e. low and high values are
simultaneously reached by both systems at the same successive iterations.

\item There is a high degree of correlation between the successive values reached by
both systems. The correlation coefficient is in most cases greater than 0.9, even if
both systems have evolved in chaotic attractors.

\item The system of less sensitivity go always ahead of the other. That is to say, the
values successively reached by the less sensitive system are always greater than those
reached by the other. This result coincides with St Matthew Theorem, a formal
conclusion derived from independent thermodynamical considerations (internal entropy
production) \cite{Leon1990}.

\item Systems evolves in an attractor that may be periodic, chaotic, fractal,
orbital or asymptotic. Low period attractors are perhaps the most frequent in NN
interactions.

\item Both systems evolves always in the same type of attractor.

\item Each point of a periodic attractor has its own independent branch in the CT.
Although periodic attractors can also result from the degeneration of a complex
attractor (chaotic, fractal or orbital).

\item When systems fall in orbital attractors, the orbit of the less sensitive system
always raises over the orbit of the more sensitive one. It is also of less size, which
means that the less sensitive system is more stable.

\item Systems can evolve from a chaotic regime to a single periodic attractor.

\item Systems can evolve from order to chaos after traversing their respective regular
and symmetrical CT. Or in other words, after a long stage of periodic behaviour the
interacting systems can evolve to a chaotic regime.

\item Systems can evolve from an initial chaotic regime to other different final
chaotic regime. Between both chaotic regimes, systems exhibit a regular periodic
behaviour while traversing their respective CT.

\item Chaotic attractors extend on a broad region and exhibit a typical internal
structure of parabolic gaps.

\item In certain cases, both systems have similar CT although they are traversed in
opposite senses (trajectory inversions).

\item The extinction of a system (to reach the value of 0) is extremely rare.

\end{enumerate}

\section{Positive-negative interactions}

\noindent According to (\ref{eqn:neg-neg 1})-(\ref{eqn:neg-neg 2}) the effects of PN
interactions are different in both interacting system: one of them benefits from the
growth of the other while this other is negatively affected by the growth of the
first. This type of interaction, therefore, models prey-predator interactions. As in
the case of NN interactions, the coevolution of both systems also depends on four real
variables within the real interval $(0, 1)$: the sensitivities and the initial seeds
of both systems. Consequently there is also here a huge number of possibilities to
explore. The behaviour diversity is now even greater than in the NN case. Some of the
most relevant characteristics of this type of recursive interaction are the following:

\begin{enumerate}

\item The coevolution trajectories of both systems are similar, although the length of
the branches may be different.

\item The coevolution trajectories can be convergent or divergent, or may suddenly
change from convergent to divergent (or viceversa). As in the NN case, they are
symmetrical with respect to the bisector $x_{n+1} = x_n$ in the case of system $X$, or
$y_{n+1} = y_n$ in the case of system $Y$.

\item The system suffering the positive interaction (predator system) is more affected
by the interaction and may evolve to extinction. In these cases the other system (prey
system) remains in a chaotic attractor.

\item It has not been observed that both systems evolves in a chaotic attractor.

\item Both systems may evolve from a chaotic regime to any other non chaotic attractor.

\item Being on the brink of extinction, the predator systems may suddenly recover and
evolve to an stable non chaotic attractor.

\item For high values of sensitivities, the branches of both CT are either radial or
spiral.

\item Spiral trajectories may be more or less dense, and the number of their
corresponding branches is variable.

\item The center of spiral trajectories is a single period-1 attractor.

\item As in the cases of natural prey-predator interactions, systems oscillate
asynchronously, so that when a system reaches a high value the other reaches a low
one, and viceversa.

\item Contrarily to the NN case, trajectory inversions have not been observed.

\item Single attractors of period-1 are frequent. They may be reached through one or
more (convergent) branches.

\item Some orbital attractors exhibit very complicated forms in their regions of
maximum curvature (folded orbits and complex orbits).

\item Orbital attractors of the systems suffering the negative interaction always
raises over the orbital attractor of the system suffering the positive interaction.
They are also of less size than it.

\end{enumerate}

\section{Positive-positive interactions}

\noindent PP interactions model pure cooperation. There are sufficient experimental
reasons to conclude that in all cases both systems evolve to extinction. The only
exception is the elemental case $S_x$ = $S_y$ = 1; $x_0$ = $y_0$ = 0.5, which evolve
to the single attractor $0.5$. Even in the case of multiple interactions systems
evolve to extinction if the number of PP interactions is sufficiently high. From
thermodynamic analysis we know, on the other hand, that this type of interaction does
not have asymmetrical coevolution trajectories of minimum entropy production. The only
trajectory of minimum entropy in this case is the bisector $y = x$.

\section{Multiple interactions}

\noindent Nature is not a set of isolated couples of systems that interact with each
other. It is rather composed of an immense and complex interaction network in which
participate millions of systems that interact with many other different systems. It
makes sense therefore to consider the possibility of analyzing more complicated
situations than the simple binary cases we have just examined. As we have seen, binary
logistic interactions can be easily generalized to any number of interacting systems
(\ref{eqn:2-interaccion m-aria}). In these conditions, attractors theorem no longer
holds and the theoretical study is much more complicated than in the binary case.
Fortunately, we can make use of computer applications to analyze the evolution of
complex interaction networks.

With the aid of the application above mentioned, some experimental work have been
performed involving up to 1000 interacting systems, being each system capable of
interacting (positively and negatively) with a variable number of other systems,
including the case of all the others. For the sake of clarity, let us term stable to
those systems that have evolved in a very low period attractor (usually a single
attractor); and network stabilizing capacity to the number $1/r$ where $r$ is the
number of iterations which are necessary for all systems to become stable. In these
conditions, the main and more remarkable conclusion we immediately get from the above
experimental analysis is the great stabilizing capacity of complex interaction
networks. By way of example, in a set of 1000 systems in which each system interacts
positively with 100 systems and negatively with other 400, all systems become stable
in less than 40 recursive logistic interactions. Next, I resume the most interesting
conclusions of the performed experiments:

\begin{enumerate}

\item The stabilizing force of the interaction network grows as the number of
interacting systems grows; and for a definite number of systems, as the number of
interactions per system grows, and as the number of PN interactions grows.

\item Collective behaviour has been observed in some interactions. In these cases
systems oscillate synchronically according to a periodic pattern.

\item The standard deviation of the network stabilizing capacity decrease as the number
of systems increases and as the number of interactions per system increases.

\item The number of stable systems increases quickly once the first system becomes
stable.

\item The final average $\overline{x}_i$ of all systems depends on the particular
experiment, but is always around a value of 0.6. This number decreases towards 0,499
as the number of positive interactions increases.

\item The standard deviation of $\overline{x}_i$ is always very low, usually less than
0.00015.

\item Complex networks of PN interactions are extremely stable. The average value of
systems in this case is 0,499, with a standard deviation less than 0.00007.

\end{enumerate}

\clearpage
\newpage

\noindent The following figures are examples of coevolution trajectories and
attractors in both NN and PN interactions. Points $(x_n, x_{n+1})$ of system X and
points $(y_n, y_{n+1})$ of system Y are plotted in red and blue respectively;
attractors are plotted in black (for more information and figures download
\cite{Leon2008})

\begin{figure}[!htb]
\centering \includegraphics[scale=1.2]{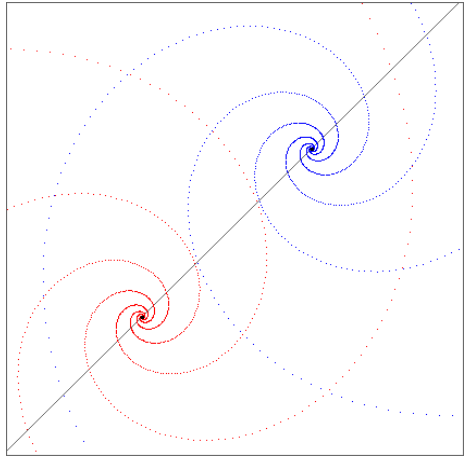} \caption{Spiral CT in a case of PN
interaction ($S_x = 0,999999; S_y$ = 0,99; $x_0$ = 0,9; $y_0$= 0,9; zoom: 28419).}
\end{figure}

\begin{figure}[!htb]
\centering \includegraphics[scale=1.2]{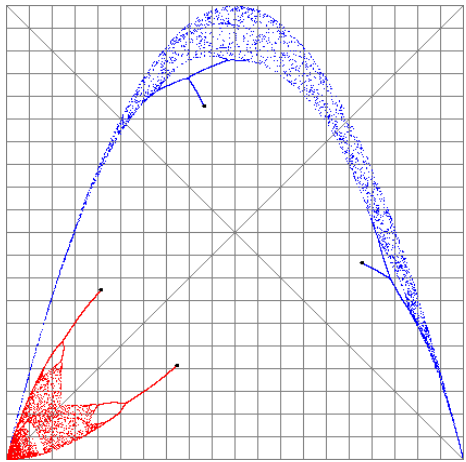} \caption{Left: From Chaos to order in a
NN interaction ($S_x$ = 0,9998; $S_y$ = 0,999; $x_0$ = 0,001; $y_0$ = 0,9). Initially
both systems evolve in a chaotic regime -parabolic regions- and then, after some
thousands interactions, they evolve through their respective symmetrical coevolution
trajectories towards a final 2-period attractor.}
\end{figure}

\begin{figure}[!htb]
\centering \includegraphics[scale=1.2]{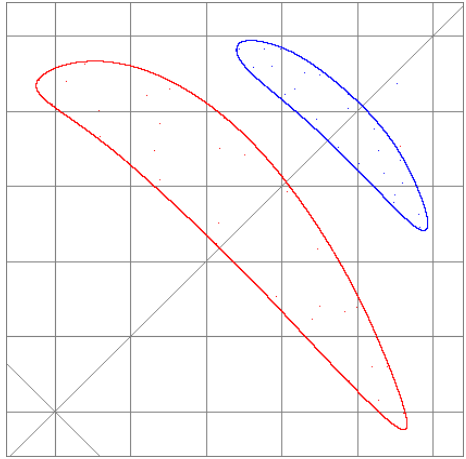} \caption{Orbital attractor of single
orbit in he interaction defined by $S_x$ = 0,988; $S_y$ = 0,3; $x_0$ = 0,9; $y_0$ =
0,9; zoom: 4.}
\end{figure}

\begin{figure}[!htb]
\centering \includegraphics[scale=1.2]{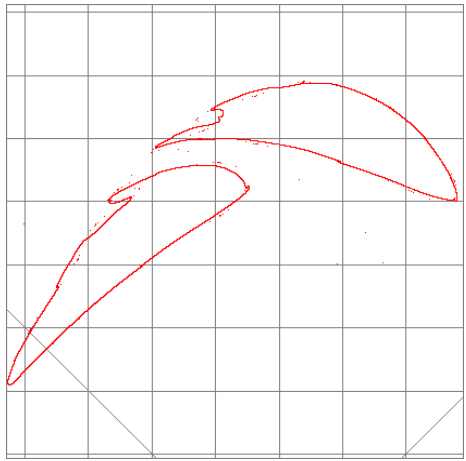} \caption{Orbital attractor of irregular
orbit. NN interaction $S_x$ = 0,335; $S_y$ = 0,333; $x_0$ = 0,3; $y_0$ = 0,2; zoom:
3.}
\end{figure}

\begin{figure}[!htb]
\centering \includegraphics[scale=1.2]{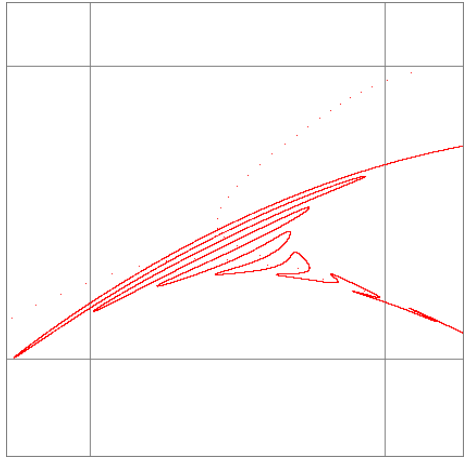} \caption{Orbital attractors of
irregular orbit. PN interaction $S_x$ = 0,9999; $S_y$ = 0,23308; $x_0$ = 0,2; $y_0$ =
0,4; zoom: 13.}
\end{figure}

\begin{figure}[!htb]
\centering \includegraphics[scale=1.2]{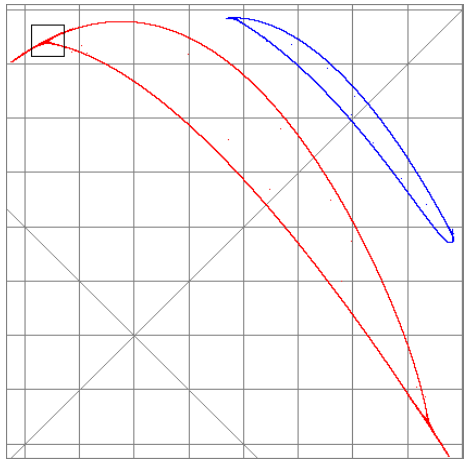} \caption{Orbital attractor of folded
orbit corresponding to the PN interaction $S_x$ = 0,999; $S_y$ = 0,232; $x_0$ = 0,9;
$y_0$ = 0,9; zoom: 3}
\end{figure}

\begin{figure}[!htb]
\centering \includegraphics[scale=1.2]{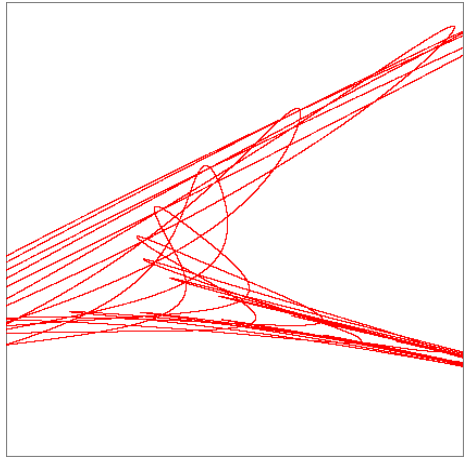} \caption{Complex region of an orbital
attractor of folded orbit corresponding to the PN interaction $S_x$ = 0,999; $S_y$ =
0,232; $x_0$ = 0,9; $y_0$ = 0,9; zoom: 97.}
\end{figure}

\begin{figure}[!htb]
\centering \includegraphics[scale=1.2]{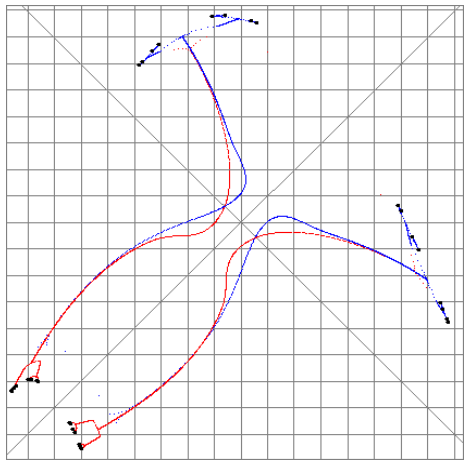} \caption{Diagonal symmetry towards a
final period-16 attractor corresponding to the NN interaction defined by $S_x$ =
0.99988; $S_y$ = 0.9976; $x_0$ = 0.765; $y_0$ = 0.234; zoom: 2.}
\end{figure}

\begin{figure}[!htb]
\centering \includegraphics[scale=1.2]{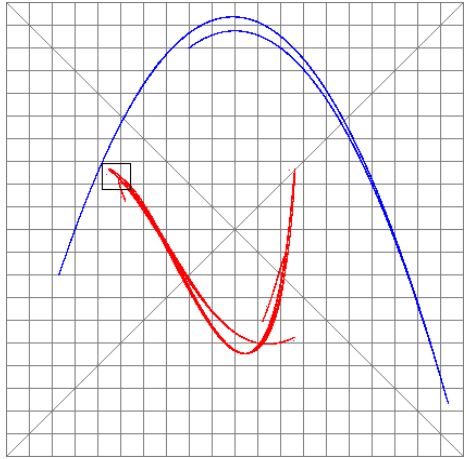} \caption{Fractal orbital corresponding
to the NN interaction $S_x$ = 0.8; $S_y$ = 0.1; $x_0$ = 0.22; $y_0$ = 0.12.}
\end{figure}

\begin{figure}[!htb]
\centering \includegraphics[scale=1.2]{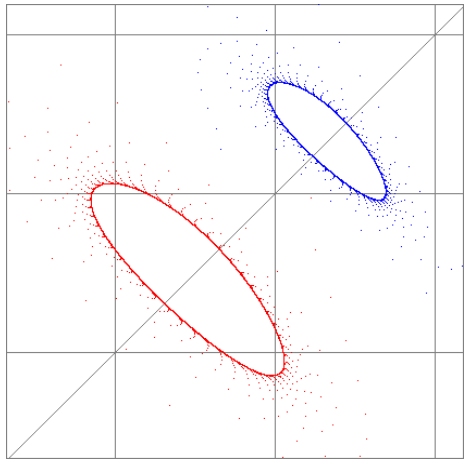} \caption{Radial CT towards an orbital
attractor of $S_x$ = 0.999; $S_y$ = 0.393999; $x_0$ = 0.45; $y_0$ = 0.67; zoom: 7}
\end{figure}

\begin{figure}[!htb]
\centering \includegraphics[scale=1.2]{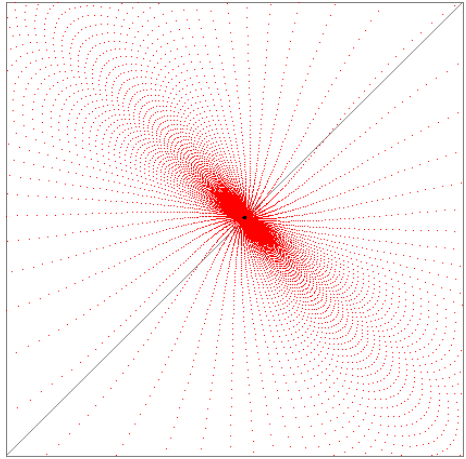} \caption{Radial CT towards a period-1
attractor of the PN interaction defined by $S_x$ = 0.999; $S_y$ = 0.399; $x_0$ = 0.8;
$y_0$ = 0.8. zoom: 115.}
\end{figure}

\clearpage
\newpage

\providecommand{\bysame}{\leavevmode\hbox to3em{\hrulefill}\thinspace}
\providecommand{\MR}{\relax\ifhmode\unskip\space\fi MR }
\providecommand{\MRhref}[2]{%
  \href{http://www.ams.org/mathscinet-getitem?mr=#1}{#2}
} \providecommand{\href}[2]{#2}

\end{document}